\documentclass[12pt, a4paper]{amsart}
 \usepackage{amsaddr}
\usepackage[onehalfspacing]{setspace}
 \usepackage{stmaryrd}
\usepackage{fullpage}
\usepackage{amsmath}
\usepackage{epigraph}
\usepackage{amsmath, amscd}
\usepackage{IEEEtrantools}
\usepackage{amssymb, tikz}
\usepackage{mathabx}
\usepackage{mathtools}
\usepackage{comment}

\calclayout

\newtheorem{thm}{Theorem}[section]

\theoremstyle{definition}

\theoremstyle{definition}

\theoremstyle{remark}

\numberwithin{equation}{section}

\def\G{\Gamma}

\newcommand{\R}{\mathbb{R}}  
\newcommand{\Z}{\mathbb{Z}}

\newcommand{\C}{\mathbb{C}}

\newcommand{\SL}{\mathrm{SL}_2 (\Z)}

\newcommand{\uh}{\mathbb{H}}

\newcommand{\sgn}{\mathrm{sign}}

\newcommand{\tr}{\mathrm{tr}}

\newcommand\SmallMatrix[4]{\ensuremath{\left(\begin{smallmatrix}#1 & #2 \\ #3 & #4\end{smallmatrix}\right)}}

\def\XXint#1#2#3{{\setbox0=\hbox{$#1{#2#3}{\displaystyle\displaystyle\int}$ }
\vcenter{\hbox{$#2#3$ }}\kern-.6\wd0}}

\usepackage{blindtext}
\usepackage{soul}

\usepackage{etoolbox}

\makeatletter
\patchcmd{\@setauthors}{\MakeUppercase\@author}{}{}{}
\makeatother

\begin{document}

\title{Cycle Integrals of the Parson Poincar\'e Series and Intersection Angles of Geodesics on Modular Curves}
\author[A. L\"ageler, M. Schwagenscheidt]{Alessandro Lägeler and Markus Schwagenscheidt}
\maketitle

\begin{abstract}
We prove a geometric formula for the cycle integrals of Parson's weight $2k$ modular integrals in terms of the intersection angles of geodesics on modular curves. Our result is an analog for modular integrals of a classical formula for the cycle integrals of certain hyperbolic Poincar\'e series, due to Katok. On the other hand, it extends a recent geometric formula of Matsusaka and Duke, Imamo\={g}lu, and T\'oth for the cycle integrals of weight 2 modular integrals. 
\end{abstract}

\section{Introduction and statement of the main results}

A classical result of Katok \cite{katok} states that for integers $k \geq 2$ the space $S_{2k}(\Gamma)$ of cusp forms of weight $2k$ for a cofinite discrete subgroup $\Gamma \subset \mathrm{SL}_2(\R)$ is generated by a family of Poincar\'e series associated with primitive hyperbolic matrices $\gamma \in \Gamma$. Explicitly, these Poincar\'e series are defined by\footnote{We use a slightly different normalization than Katok \cite{katok} to simplify our formulas.}
\begin{align}\label{fkgamma}
f_{k,\gamma}(z) = - \frac{D_\gamma^{k-\frac{1}{2}}}{\pi}\sum_{g \in  \Gamma_\gamma \backslash\Gamma}\frac{1}{(Q_{\gamma }\circ g)(z,1)^{k}},
\end{align}
where $\Gamma_\gamma = \{\pm \gamma^n: n \in \Z\}$\footnote{We will assume throughout that $-1 \in \Gamma$.}, $Q_\gamma(x,y) = cx^2+(d-a)xy-by^2$ is the binary quadratic form corresponding to $\gamma = \left(\begin{smallmatrix}a & b \\ c & d \end{smallmatrix}\right)$, $D_\gamma = \tr(\gamma)^2-4$ is the discriminant of $Q_\gamma$, and $\Gamma$ acts on binary quadratic forms in the usual way. These hyperbolic Poincar\'e series have many other interesting applications, the most prominent one being Kohnen's \cite{kohnen} construction of the holomorphic kernel function for the Shimura correspondence.

Katok \cite{katok} also gave a beautiful geometric formula for (the imaginary part of) the geodesic cycle integrals of the cusp forms $f_{k,\gamma}(z)$. If $\gamma,\sigma \in \Gamma$ are primitive hyperbolic elements in $\Gamma$ which have positive trace and which are not conjugacy equivalent, and $S_{\gamma}$ denotes the geodesic semi-circle in $\uh$ connecting the two real fixed points of $\gamma$, then Katok's formula \cite[Theorem~3]{katok} states that 
\begin{align} \label{katokformula}
    \mathrm{Im} \left(\int_{z_0}^{\sigma. z_0} f_{k,\gamma}(z) Q_\sigma(z, 1)^{k - 1} dz \right) = (D_\gamma D_\sigma)^{\frac{k - 1}{2}} \sum_{p \in [S_\gamma]\cap [S_\sigma]} \mu_p^{k} P_{k - 1}( \cos \theta_p),
\end{align}
where $z_0 \in \uh$ and the path of integration can be chosen arbitrarily, the sum runs over the finitely many intersection points $p$ of the closed geodesics $[S_\gamma] = \Gamma_\gamma \backslash S_\gamma$ and $[S_\sigma] = \Gamma_\sigma \backslash S_{\sigma}$ in $\Gamma \setminus \uh$, and $P_r$ denotes the $r$-th Legendre polynomial. Moreover, $\theta_p = \theta_p(\gamma,\sigma) \in [0, \pi]$ denotes the \emph{intersection angle} at $p$, which is measured counterclockwise from the tangent at $S_{\gamma}$ to the tangent at $S_{\sigma}$ at $p$, and $\mu_p = \mu_p(\gamma,\sigma) \in \{\pm 1\}$ denotes the \emph{sign of the intersection} at $p$, which is defined as follows: let $g \in \Gamma$ be chosen such that the intersection point $p$ corresponds to the intersection point of $S_\gamma$ and $S_{g \sigma g^{-1}}$ in $\uh$, and suppose that $S_\gamma$ and $S_{g \sigma g^{-1}}$ are oriented clockwise (which means that the lower left entries of $\gamma$ and $g\sigma g^{-1}$ are positive). Then $\mu_p(\gamma,\sigma) = +1$ if the left endpoint of $S_{g \sigma g^{-1}}$ lies between the two endpoints of $S_\gamma$, and  $\mu_p(\gamma,\sigma) = -1$ otherwise. The sign of $\mu_p$ changes if the orientation of either $S_\gamma$ or $S_{g \sigma g^{-1}}$ is reversed. Note that $\theta_p(\gamma,\sigma)$ does not depend on the orientation of $S_\gamma$ or $S_\sigma$, but it does depend on the order of $\gamma,\sigma$, that is, we have $\theta_p(\sigma,\gamma) = \pi - \theta_p(\gamma,\sigma)$. Similarly, we have $\mu_p(\gamma,\sigma) = -\mu_p(\sigma,\gamma)$.

More recently, Matsusaka~\cite{matsusaka} investigated the (homogenized) cycle integrals of certain modular integrals of weight $2$ for $\SL$ with rational period functions. These modular integrals were constructed by Duke, Imamo\={g}lu, and T\'oth in \cite{dit,ditmodularintegrals}, and are defined\footnote{Again, our normalization differs from \cite{matsusaka,dit,ditmodularintegrals}.} for primitive hyperbolic $\gamma \in \SL$ by
\begin{align}\label{ditmodularintegral}
F_{\gamma}(z) = \frac{2\sqrt{D_\gamma}}{\pi}\sum_{n = 0}^\infty \left(\int_{z_0}^{\gamma z_0}j_n(\tau)\frac{d\tau}{Q_\gamma(\tau,1)}\right)e^{2\pi i n z},
\end{align}
where $j_n(\tau)$ denotes the unique weakly holomorphic modular function for $\SL$ whose Fourier expansion has the shape $j_n(\tau) = q^{-n}+O(q)$ with $q = e^{2 \pi i \tau}$. It was shown in \cite{dit} that the series defining $F_\gamma(z)$ converges to a holomorphic function on $\uh$, and satisfies for any $\sigma \in \SL$ the transformation formula
\begin{align}\label{ditcocycle}
r_\gamma(\sigma,z) = (F_{\gamma}|_{2}\sigma)(z) - F_{\gamma}(z) = \frac{2\sqrt{D_\gamma}}{\pi}\sum_{\substack{g \in \Gamma_\gamma\backslash \Gamma \\ w_{g^{-1}\gamma g}' < \sigma^{-1}. i \infty < w_{g^{-1}\gamma g}}}\frac{\sgn(Q_\gamma\circ g)}{(Q_\gamma\circ g)(z,1)},
\end{align}
where $w_\gamma' < w_\gamma$ denote the two real fixed points of $\gamma$, and we put $\sgn(Q) = \sgn(A)$ for a binary quadratic form $Q(x, y) = Ax^2 + Bxy + Cy^2$. Note that the sum on the right-hand side is finite. The function $\sigma \mapsto r_\gamma(\sigma,z)$ defines a holomorphic weight $2$ cocycle for $\SL$ with values in the rational functions on $\C$, and the function $F_\gamma(z)$ is called a \emph{modular integral} of weight $2$ for $r_{\gamma}(\sigma,z)$. Matsusaka proved the remarkable formula
\begin{align}\label{matsusakaformula}
    \mathrm{Im}\left(\lim_{n \to \infty}\int_{\sigma^n. z_0}^{\sigma^{n+1}.z_0}F_{\gamma}(z)dz\right) = \sum_{p \in [S_\gamma]\cap[S_\sigma]}1.
\end{align}
where $z_0 \in \uh$ and the path of integration can be chosen arbitrarily, and $\gamma,\sigma \in \SL$ are primitive hyperbolic matrices with positive trace which are not conjugacy equivalent (compare \cite[Corollary~3.8, Theorem~3.3]{matsusaka}). Note that the homogenization of the cycle integral on the left-hand side is necessary to make it independent of $z_0$, and a conjugacy class invariant in $\sigma$. The right-hand side of the formula \eqref{matsusakaformula} counts the number of intersections of $[S_\gamma]$ and $[S_\sigma]$ in $\SL \backslash \uh$, which by \cite[Theorem~3]{ditlinking} can also be interpreted as the linking number of certain modular knots associated with $\gamma$ and $\sigma$.

Notice that there is a striking similarity between Matsusaka's formula \eqref{matsusakaformula} and (the formal specialization to $k = 1$ of) Katok's formula \eqref{katokformula}.
Motivated by this observation, in the present work we extend Matsusaka's formula \eqref{matsusakaformula} to certain modular integrals of higher weight $2k$ (with $k \geq 2$) for cofinite discrete subgroups $\Gamma \subset \mathrm{SL}_2(\R)$, by evaluating their homogenized cycle integrals in terms of intersection angles of geodesics in $\Gamma \backslash \uh$, much in the spirit of Katok's formula \eqref{katokformula}.

The modular integrals we consider here were introduced by Parson\footnote{In fact, our definition is not precisely Parson's, but differs from her original Poincar\'e series by a cusp form. Moreover, we use a different normalization than Parson to match our normalization of \eqref{fkgamma}.} in \cite{parson}, and are defined for integers $k \geq 2$ and primitive hyperbolic $\gamma \in \Gamma$ by
\begin{align}\label{parsonmodularintegral}
F_{k,\gamma}(z) := - \frac{D_\gamma^{k-\frac{1}{2}}}{\pi}\sum_{g \in \Gamma_\gamma \backslash\Gamma}\frac{\sgn(Q_\gamma \circ g)}{(Q_{\gamma }\circ g)(z,1)^{k}}.
\end{align}
A direct computation shows that the holomorphic function $F_{k,\gamma}(z)$ satsifies for any $\sigma\in \Gamma$ the transformation law
\begin{align}\label{parsoncocycle}
r_{k,\gamma}(\sigma,z): = (F_{k,\gamma}|_{2k}\sigma)(z) - F_{k,\gamma}(z) = \frac{2D_\gamma^{k-\frac{1}{2}}}{\pi}\sum_{\substack{g \in \Gamma_\gamma \backslash\Gamma \\ w_{g^{-1}\gamma g}' < \sigma^{-1}. i \infty < w_{g^{-1}\gamma g}}}\frac{\sgn(Q_\gamma\circ g)}{(Q_\gamma\circ g)(z,1)^{k}}.
\end{align}
The sum on the right hand side is finite. In particular, the map $\sigma \mapsto r_{k,\gamma}(\sigma,z)$ is a holomorphic weight $2k$ cocycle for $\Gamma$ with values in the rational functions on $\C$, and $F_{k,\gamma}(z)$ is a modular integral for $r_{k,\gamma}(\sigma,z)$. 

Notice that the cocycle $r_{\gamma}(\sigma,z)$ in \eqref{ditcocycle} is the specialization to $k = 1$ of the cocycle $r_{k,\gamma}(\sigma,z)$ in \eqref{parsoncocycle}. Hence, we may view the weight $2$ modular integral $F_{\gamma}(z)$ defined in \eqref{ditmodularintegral} as the $k = 1$ analog of Parson's modular integral $F_{k,\gamma}(z)$ defined in \eqref{parsonmodularintegral}. However, $F_{k,\gamma}(z)$ does not converge for $k = 1$, although it is probably possible to extend the definition \eqref{parsonmodularintegral} to $k = 1$ using Hecke's trick, and to show that $F_{1,\gamma}(z) = F_{\gamma}(z)$.

Our main result is the following geometric formula for (the imaginary part of) the cycle integrals of Parson's modular integrals $F_{k,\gamma}(z)$. It is an analog of Katok's formula \eqref{katokformula} for modular integrals, and a higher weight analog of Matsusaka's formula \eqref{matsusakaformula}.

\begin{thm}\label{mainresult}
Let $k \geq 2$ be an integer and let $\gamma$, $\sigma$ be primitive hyperbolic elements in $\Gamma$ with positive trace which are not conjugacy equivalent. We have 
\begin{equation} \label{modintkatok}
    \mathrm{Im} \left( \lim_{n \to \infty}\int_{\sigma^n.z_0}^{\sigma^{n+1}.z_0} F_{k,\gamma}(z) Q_\sigma(z, 1)^{k - 1} dz \right) = (D_\gamma D_\sigma)^{\frac{k - 1}{2}} \sum_{p \in [S_\gamma]\cap [S_\sigma]} \mu_p^{k - 1} P_{k - 1}( \cos \theta_p),
\end{equation}
where the notation is as in \eqref{katokformula}.
\end{thm}

The proof of Theorem~\ref{mainresult} will be given in Section~\ref{proofs} below. We would like to give a quick proof of the fact that the limit on the left-hand side exists, and is independent of $z_0$. First, note that the limit of $\sigma^n.z_0$ as $n \to \infty$ is independent of the choice of $z_0 \in \uh$ and converges to $w_\sigma$ (if we assume for the moment that $\sigma$ has positive trace and positive lower left entry; see \cite[Lemma 2.7]{matsusaka}). Now the function $x\mapsto\displaystyle\int_{x}^{\sigma.x} F_{k,\gamma}(z) Q_\sigma(z, 1)^{k - 1} dz$ is Lipschitz-continuous when $x$ is approaching $w_\sigma$. Indeed, note that \begin{align} \begin{split}\label{Lipschitz}
    &\left\vert \int_{x_0}^{\sigma.x_0} F_{k,\gamma}(z) Q_\sigma(z, 1)^{k - 1} dz - \int_{x_1}^{\sigma.x_1} F_{k,\gamma}(z) Q_\sigma(z, 1)^{k - 1} dz \right\vert \\
    &= \left\vert \int_{x_0}^{x_1} r_{k,\gamma} (\sigma, z) Q_\sigma(z, 1)^{k - 1} \, dz \right\vert \leq L_{k, \gamma, \sigma} \vert x_0 - x_1 \vert,
\end{split}
\end{align}
for some constant $L_{k,\gamma,\sigma} > 0$ and $x_0,x_1$ close to $w_\sigma$, as $r_{k,\gamma}(\sigma, z)$ is holomorphic at $w_\sigma$ if $\gamma$ and $\sigma$ are not conjugacy equivalent, so $\vert r_{k,\gamma}(\sigma, z) Q_\sigma(z, 1)^{k - 1} \vert$ is bounded in a neighbourhood of $w_\sigma$. This implies that $\displaystyle\int_{\sigma^n.z_0}^{\sigma^{n+1}.z_0} F_{k,\gamma}(z) Q_\sigma(z, 1)^{k - 1} dz$ is a Cauchy sequence. Moreover, if we put $x_0 = \sigma^n z_0$ and $x_1 = \sigma^n z_1$ in \eqref{Lipschitz} and take the limit as $n \to \infty$, we see that the left-hand side in Theorem~\ref{mainresult} is independent of $z_0$. Note that this also implies that the homogenized cycle integral of $F_{k,\gamma}(z)$ is a conjugacy class invariant in $\sigma$. We would also like to remark that the right-hand side in Theorem~\ref{mainresult} is a finite sum, which can be explicitly computed (numerically) as explained in \cite{rickards}.

%


Previous to Matsusaka's formula (\ref{matsusakaformula}), Duke, Imamo\={g}lu, and T\'oth proved its "parabolic version" which expresses the number of intersections of the net of the geodesics equivalent to $S_\gamma$ with the non-compact geodesic $S_{- d / c}$ from $- \frac{d}{c}$ to $i \infty$ as the central value of the twisted $L$-function 
\begin{align*}
L_{\gamma}\left( s, - \frac{d}{c} \right) = \sum_{n = 1}^\infty \frac{a_\gamma(n) e^{ - 2\pi i\frac{d}{c} n}}{n^s}, \quad (c,d) = 1, \; c > 0, \quad \mathrm{Re}(s) \gg 1,
\end{align*}
of the modular integral $F_\gamma(z) = \sum_{n = 1}^\infty a_\gamma(n) e^{2 \pi i n z}$ in (\ref{ditmodularintegral}). The explicit formula reads \begin{equation} \label{ditintersectionformula}
   \frac{1}{2\pi} \; \mathrm{Re} L_{\gamma}\left(1, - \frac{d}{c} \right) = \sum_{p \in [S_\gamma] \; \cap \; S_{- d / c}} 1,
\end{equation}
where the equivalence class is over the action of $\SL$ on the geodesic $S_\gamma$ (see \cite[Theorem 5.2]{ditlinking}). We generalize (\ref{ditintersectionformula}) to the higher weight case. 

The Parson Poincar\'e series $F_{k,\gamma}(z)$ defined in \eqref{parsonmodularintegral} has a Fourier expansion of the shape
\[
F_{k,\gamma}(z) = \sum_{\substack{n \in \frac{1}{N}\Z, n > 0}}a_{k,\gamma}(n)e^{2\pi i n z},
\]
where $N$ denotes the width of the cusp $\infty$ with respect to the group $\Gamma$. The Fourier coefficients $a_{k,\gamma}(n)$ can be explicitly computed in terms of Kloosterman sums and Bessel functions as in \cite[Theorem~3.1]{parson}, or in terms of cycle integrals of weakly holomorphic modular forms as in \cite[Theorem~3]{ditmodularintegrals}. We define the twisted $L$-function of $F_{k,\gamma}(z)$ by
\begin{align*}
L_{k,\gamma}\left( s, - \frac{d}{c} \right) = \sum_{n \in \frac{1}{N}\Z, n > 0}\frac{a_{k,\gamma}(n) e^{- 2\pi i\frac{d}{c} n}}{n^s}, \quad (c,d) = 1, \; c > 0, \quad \mathrm{Re}(s) \gg 1,
\end{align*}
At its central value, the twisted $L$-function satisfies the following analog of (\ref{ditintersectionformula}).

\begin{thm}\label{mainresult2}
Let $k \geq 2$ be an odd integer, $\gamma$ be a hyperbolic element with positive trace and $d, c$ be comprime integers with $c > 0$. We have \begin{align} \label{intersectionformulahigherk}
    (-1)^{\frac{k - 1}{2}} \frac{(k - 1)!}{(2\pi)^k} \; \mathrm{Re} L_{k,\gamma}\left(k, - \frac{d}{c} \right) = D_\gamma^{\frac{k - 1}{2}} \sum_{p \in [S_\gamma] \; \cap \;  S_{- d / c}} P_{k - 1} \left( \cos \theta_p \right),
\end{align}
where $\theta_p$ denotes the angle of intersection at the point $p$ between the geodesic in the equivalence class of $S_\gamma$ and the geodesic $S_{- d / c}$. 
\end{thm}

The left-hand side of \eqref{intersectionformulahigherk} can be written as the imaginary part of the cycle integral of $F_{k,\gamma}(z)$ along the non-compact geodesic from $-d/c$ to $i\infty$, so Theorem~\ref{mainresult2} can be viewed as a "parabolic" analog of Theorem~\ref{mainresult}.


\section{The proofs of Theorem~\ref{mainresult} and Theorem~\ref{mainresult2}}\label{proofs}

\subsection{Proof of Theorem~\ref{mainresult}}
    The proof is similar to the proof of Katok's formula \eqref{katokformula}, compare \cite[Theorem~3]{katok}. The main difference is that we can't unfold the Parson Poincar\'e series as such.

    Up to interchanging $\gamma$ with $\gamma^{-1}$ and $\sigma$ with $\sigma^{-1}$ we may assume that the lower left entries of $\gamma$, $\sigma$ are positive. First, we rewrite 
    \[
    \int_{\sigma^n.z_0}^{\sigma^{n + 1}.z_0} F_{k,\gamma}(z) Q_\sigma(z, 1)^{k - 1} dz = \int_{z_0}^{\sigma.z_0} (c_n z + d_n)^{-2k} F_{k,\gamma}(\sigma^n.z) Q_\sigma(z, 1)^{k - 1} dz,
    \]
where we put $\sigma^n = \SmallMatrix{*}{*}{c_n}{d_n}$. A direct computation shows that
\begin{equation} \label{lemmaintersect}
    \sgn(Q_\gamma \circ \sigma^{-n} (1, 0)) \to \sgn(Q_\gamma(w_\sigma', 1)) \; \; \mathrm{as} \; n \to + \infty;
\end{equation}
see \cite[Lemma 2.7]{matsusaka}. If the two geodesics $S_\gamma$ and $S_\sigma$ intersect, we have $\sgn(Q_\gamma(w_\sigma', 1)) = - \mu_p(\gamma, \sigma)$.
Now the Parson Poincar\'e series becomes \begin{align*}
    (c_n z + d_n)^{-2k} F_{k,\gamma}(\sigma^n.z) &= - \frac{D_\gamma^{k - 1 / 2}}{\pi} \sum_{g \in \G_\gamma \setminus \G} \frac{\sgn(Q_\gamma \circ g)}{(c_n z + d_n)^{2k} (Q_\gamma \circ g)(\sigma^n.z, 1)^k} \\
    &= - \frac{D_\gamma^{k - 1 / 2}}{\pi}\sum_{g \in \G_\gamma \setminus \G} \frac{\sgn(Q_\gamma \circ g \sigma^{-n}))}{(Q_\gamma \circ g)(z, 1)^k} \\
    &\to - \frac{D_\gamma^{k - 1 / 2}}{\pi} \sum_{g \in \G_\gamma \setminus \G} \frac{\sgn((Q_\gamma \circ g)(w_\sigma', 1))}{(Q_\gamma \circ g)(z, 1)^k}, \; \mathrm{as} \; n \to + \infty,
\end{align*}
where in the last line we applied (\ref{lemmaintersect}). The resulting series is modular of weight $2k$ for the group $\Gamma_\sigma=\{\pm \sigma^n: n\in \Z\}$. Hence, the cycle integral $$\int_{z_0}^{\sigma.z_0} \sum_{g \in \G_\gamma \setminus \G} \frac{\sgn((Q_\gamma \circ g)(w_\sigma', 1))}{(Q_\gamma \circ g)(z, 1)^k} Q_\sigma(z, 1)^{k - 1} dz$$
is independent of the choice of $z_0$.

Now a typical unfolding argument yields \begin{align*}
    &\int_{z_0}^{\sigma.z_0} \sum_{g \in \G_\gamma \setminus \G} \frac{\sgn((Q_\gamma \circ g)(w_\sigma', 1))}{(Q_\gamma \circ g)(z, 1)^k} Q_\sigma(z, 1)^{k - 1} dz \\
    &= \sum_{g \in \G_\gamma \setminus \G / \G_\sigma} \sum_{m \in \Z} \int_{\sigma^m.z_0}^{\sigma^{m + 1}.z_0} \frac{\sgn((Q_\gamma \circ g)(w_\sigma', 1))}{(Q_\gamma \circ g) (z, 1)^k} Q_\sigma(z, 1)^{k - 1} dz \\
    &= \sum_{g \in \G_\gamma \setminus \G / \G_\sigma} \int_{S_\sigma} \frac{\sgn((Q_\gamma \circ g)(w_\sigma', 1))}{(Q_\gamma \circ g) (z, 1)^k} Q_\sigma(z, 1)^{k - 1} dz.
\end{align*}

Subtracting the complex conjugate from the last equation, gives the sum of integrals  $$\sum_{g \in \G_\gamma \setminus \G / \G_\sigma} \displaystyle\int_{C(\sigma)} \frac{\sgn((Q_\gamma \circ g)(w_\sigma', 1))}{(Q_\gamma \circ g) (z, 1)^k} Q_\sigma(z, 1)^{k - 1} dz$$ over the circle $C(\sigma)$ through the roots of $Q_\sigma(z, 1) = 0$. The integrands are all meromorphic, with poles only at the real roots of $(Q_\gamma \circ g)(z,1)$. Hence, by Cauchy's theorem, they vanish if the geodesic connecting the roots of $(Q_\gamma \circ g) (z,1)$ does not intersect $S_\sigma$. Thus, we are left with
\[
\frac{D_\gamma^{k - 1 / 2}}{\pi} \sum_{p \in [S_\gamma] \cap [S_\sigma]} \mu_p(\gamma, \sigma) \displaystyle\int_{C(\sigma)} \frac{1}{(Q \circ g) (z, 1)^k} Q_\sigma(z, 1)^{k - 1} dz,
\]
as $\mu_p(\gamma, \sigma) = - \sgn((Q \circ g) (w_\sigma', 1))$ if the geodesics $S_\sigma$ and $S_{g^{-1}\gamma  g}$ intersect. 

The integrals evaluate to $$\displaystyle\int_{C(\sigma)} \frac{1}{(Q_\gamma \circ g) (z, 1)^k} Q_\sigma(z, 1)^{k - 1} dz = 2 \pi i D_\gamma^{- k / 2} D_\sigma^{\frac{k - 1}{2}} \mu_p^k P_{k - 1} (\cos \theta_p);$$
see \cite[p. 478]{katok}. This finishes the proof of Theorem~\ref{mainresult}.

\subsection{Proof of Theorem~\ref{mainresult2}}
The proof is a careful application of \cite[Lemma 2]{katok}, which asserts that for odd $k \in \Z, k \geq 3$, and any $A, B, C \in \R$ with $D = B^2 - 4AC > 0$, we have \begin{equation} \label{legendreeva}
    \int_{-\infty}^{\infty} \frac{t^{k - 1}}{(-At^2 + Bit + C)^k} dt = \begin{cases}0, &AC > 0, \\ (-1)^{\frac{k + 1}{2}} \sgn(A) 2 \pi D^{- k / 2}  P_{k - 1} \left( \frac{B}{\sqrt{D}} \right), &AC < 0. \end{cases}
\end{equation}

Consider the integral $\displaystyle\int_{- \frac{d}{c}}^{i \infty} F_{k, \gamma}(z) (cz + d)^{k - 1} dz$. A standard argument gives that \begin{equation*}
    \int_{- \frac{d}{c}}^{i \infty} F_{k, \gamma}(z) (cz + d)^{k - 1} dz = \left( \frac{c}{2 \pi } \right)^k \Gamma(k) \frac{i^k}{c} L_{\gamma, k}\left( k, -\frac{d}{c} \right).
\end{equation*} 
The imaginary part of the integral is equal to \begin{align*}
    \mathrm{Im} \int_{- \frac{d}{c}}^{i \infty} F_{k, \gamma}(z) (cz + d)^{k - 1} dz &= c^{k - 1} \mathrm{Im}  \int_{0}^{i \infty} F_{k, \gamma}\left(z - \frac{d}{c} \right) z^{k - 1} dz \\
    &= - c^{k - 1} \frac{D_\gamma^{k-\frac{1}{2}}}{\pi} \sum_{g \in \G_\gamma \setminus \G} \mathrm{Im} \left( \int_{0}^{i \infty} \frac{\sgn(Q_\gamma \circ g) z^{k - 1}}{(Q_\gamma \circ g)(z - d / c, 1)^k} dz \right),
\end{align*}
where we can exchange sum and integral by Fubini's theorem, as the integrand is absolutely convergent.

Fix $g \in \G_\gamma \setminus \G$ for the moment and write $Az^2 + Bz + C$ for $(Q_\gamma \circ g)(z, 1)$. The quadratic form 
\[
(Q_\gamma\circ g)(z - d / c, 1) = Az^2 + \left(B - 2 A d / c \right)z + \left( A \left( d / c \right)^2 - B d / c + C \right) = A'z^2 + B'z + C'
\]
intersects with the non-compact geodesic $S_{- d / c}$ if and only if 
\[
A' C' = A\left( A \left( d / c \right)^2 - B d / c + C \right) < 0.
\]
With (\ref{legendreeva}), we get 
\begin{align*}
    \mathrm{Im} \left( \int_{0}^{i\infty} \frac{\sgn(A') z^{k - 1}}{\left(A'z^2 + B'z + C'\right)^k} dz \right) 
    &= \frac{(-1)^{\frac{k - 1}{2}}}{2} \int_{-\infty}^{\infty} \frac{\sgn(A') t^{k - 1}}{\left(-A't^2 + B'it + C'\right)^k} dt  \\
    &= - \pi D_\gamma^{- k / 2}  P_{k - 1} \left( \frac{B'}{\sqrt{D_\gamma}} \right)
\end{align*}
if and only if the geodesic $S_{g\gamma g^{-1}}$ intersects the non-compact geodesic $S_{- d / c}$ associated to the quadratic form $cz + d$. Otherwise, the integral evaluates to zero by (\ref{legendreeva}). By \cite[Proposition~2.2]{rickards}, the intersection angle $\theta_p$ between these two geodesics is given by $\cos \theta_p = \frac{Bc - 2Ad}{c \sqrt{D_\gamma}}$. This finishes the proof.

\section{Additional Remarks}

In this section, we present some further properties and possible applications of the homogenized cycle integrals and periods of Parson's modular integrals $F_{k,\gamma}(z)$.

\subsection{Explicit representation of the cycle integral} As the Parson Poincar\'e series is no longer modular, the cycle integral $\displaystyle\int_{z_0}^{\sigma.z_0} F_{k, \gamma} (z) Q_\sigma(z, 1)^{k - 1} dz$ depends on the choice of the point $z_0$ and is a complicated function in terms of $z_0$. Only when taking the homogenization we get a conjugacy class invariant object independent of the choice of $z_0$, which has the nice representation on the right hand side of (\ref{modintkatok}).

Explicitly, we have \begin{align*}
    &\int_{z_0}^{\sigma.z_0} F_{k, \gamma}(z) Q_\sigma(z, 1)^{k - 1} dz = \int_{i\infty}^{\sigma.i \infty} F_{k, \gamma}(z) Q_\sigma(z, 1)^{k - 1} dz \\
    &+ \sum_{n = 0}^{2k - 2} \sum_{\substack{g \in \G_\gamma \setminus \G, \\ w_{g \gamma g^{-1}}' < \sigma^{-1}.i \infty < w_{g \gamma g^{-1}}}} \rho_{n, g, \gamma, \sigma}(z_0) \;_2 F_1\left(k, 2k - 1 - n; 2k; 1 - \frac{z_0 - w_Q}{z_0 - w_Q'} \right),
\end{align*}
where $\rho_{n, g, \gamma, \sigma}(z_0)$ is a rational function given by $\rho_{n, g, \gamma, \sigma}(z_0) = \frac{\Gamma(2k - n - 1)}{\Gamma(2k)} \frac{\partial_z^n Q_\sigma(z, 1)^{k - 1} \vert_{z = z_0}}{(z_0 - w_{g \gamma g^{-1}}')^{n - 2k - 1}}$. One can also see from this representation that the integral converges as $z_0 \to w_\sigma$. 

To prove this, we use $$\int_{z_0}^{\sigma.z_0} F_{k,\gamma}(z) Q_\sigma(z, 1)^{k - 1} dz = \int_{i \infty}^{\sigma.i \infty} F_{k,\gamma}(z) Q_\sigma(z, 1)^{k - 1} dz + \int_{i \infty}^{z_0} r_{k,\gamma}(\sigma, z) Q_\sigma(z, 1)^{k - 1} dz.$$
This follows from differentiating both sides in $z_0$ and observing that both sides are equal at $z_0 \to i \infty$. Rewriting the polynomial $Q_\sigma(z, 1)^{k - 1} = \sum_{n = 1}^{2k - 2} a_{n, \sigma}(z_0) (z - z_0)^n$ in its Taylor expansion about $z_0$, we obtain a finite sum of integrals $$\int_{i \infty}^{z_0} r_{k,\gamma}(\sigma, z) Q_\sigma(z, 1)^{k - 1} dz = \sum_{n = 1}^{2k - 2} a_{n, \sigma}(z_0) \sum_{\substack{Q \sim Q_\gamma, \\ w_Q' < \sigma^{-1}.i \infty < w_Q}} \int_{z_0}^{i \infty} \frac{(z - z_0)^n}{(z - w_Q)^k(z - w_Q')^k} dz,$$ which can be solved. Standard integral transformations give \begin{align}
        \begin{split} \label{integralevaluation}
        &\int_{z_0}^{i\infty}\frac{(z-z_0)^{n}}{(z-w_Q)^k(z-w_Q')^k}dz = \\
        &\frac{\Gamma(2k - n - 1) \Gamma(n + 1)}{\Gamma(2k)} (z-w_Q')^{n - 2k + 1} \ _2 F_1\left(k, 2k - n - 1, 2k; 1-\frac{z-w_Q}{z-w_Q'}\right).
            \end{split}
\end{align}

That the right hand side of (\ref{integralevaluation}) is symmetric in $w_Q'$ and $w_Q$ also follows from the identity $\ _2F_1(c-a,b,c;z/(z-1)) = (1-z)^b \ _2 F_1(a,b,c;z)$.

The integral evaluation (\ref{integralevaluation}) can also be used to prove that the weight $2 - 2k$ cocycle \begin{align*}
    R_\gamma(\sigma, z) &= \frac{(-2 \pi i)^{2k - 1} D^{k - 1 / 2}}{\pi (2k - 1)! \binom{2k - 2}{k - 1}} \sum_{w_Q' < - \frac{d}{c} < w_Q} \frac{1}{\vert Q(1, 0) \vert (z - w_Q')} \ _2 F_1\bigg(k,1,2k;1-\frac{z-w_Q}{z-w_Q'}\bigg) \\
&+ \frac{(-2 \pi i)^{2k - 1}}{ (2k - 2)!} \frac{i}{c^{2k-1}} \sum_{n = 0}^{2k-2}\binom{2k-2}{n}i^n \left(\frac{c}{2 \pi} \right)^{n + 1}\Gamma(n + 1) L_{k, \gamma}(n+1,a/c) (cz+d)^{n}
\end{align*}
for $\sigma = \left( \begin{smallmatrix} a & b \\ c & d \end{smallmatrix} \right) \in \G$ and $Q$ running over the equivalence class $[Q_\gamma]$ is a $(2k - 1)$-th primitive of $r_\gamma(g, z)$, i.e. $\mathcal{D}^{2k - 1} R_\gamma(g, z) = r_\gamma(g, z)$ with $\mathcal{D} = \frac{1}{2\pi i}\frac{\partial}{\partial z}$. The Fourier coefficients $a_{k, \gamma}(n)$ can be explicitly calculated as in \cite[Theorem~3]{parson}.

\subsection{Periods of modular integrals} In this subsection we let $\Gamma = \SL$. Kohnen and Zagier~\cite{kohnenzagier} studied the periods of the cusp forms $f_{k,\gamma}(z)$ defined in~\eqref{fkgamma}, and showed that certain linear combinations of these periods are rational numbers.

A natural follow-up question to Theorem \ref{mainresult} would be to study the periods of the Parson Poincar\'e series, i.e. $p_n(F_{k, \gamma}) = \displaystyle\int_0^\infty F_{k, \gamma}(it) t^n \; dt$ for $0 \leq n \leq 2k-2$. As in~\cite{kohnenzagier} we define the symmetrizations
      \[
      F_{k, \gamma}^{+} = F_{k,\gamma} + F_{k, \gamma'}(z), \qquad   F_{k, \gamma}^{-} = i( F_{k, \gamma}(z) - F_{k, \gamma'}(z)),
      \]
where $\gamma' = \left( \begin{smallmatrix} -1 & 0 \\ 0 & 1 \end{smallmatrix} \right) \gamma \left( \begin{smallmatrix} -1 & 0 \\ 0 & 1 \end{smallmatrix} \right)$. We split the period polynomial $$p(F_{k, \gamma})(z) = \displaystyle\int_{0}^{i\infty} F_{k, \gamma} (z)(x-z)^{2k-2}dz = \sum_{n=0}^{2k-2} i^{-n+1}\binom{2k-2}{n}p_n(F_{k, \gamma})x^{2k-2-n}$$ of $F_{k, \gamma}$ into its even and odd part as $p(F_{k, \gamma}) = ip^+(F_{k, \gamma})+p^-(F_{k, \gamma})$ with
      \begin{align*}
          p^+(F_{k, \gamma})(x) &= \sum_{\substack{0 \leq n \leq 2k-2 \\ n \text{ even}}}(-1)^{n/2}\binom{2k-2}{n}p_n(F_{k, \gamma})x^{2k-2-n}, \\
          p^-(F_{k, \gamma})(x) &= \sum_{\substack{0 < n < 2k-2 \\ n \text{ odd}}}(-1)^{(n-1)/2}\binom{2k-2}{n}p_n(F_{k, \gamma})x^{2k-2-n}.
      \end{align*}

By closely following the proof of \cite[Theorem 5]{kohnenzagier}, one obtains
        \begin{align} \label{thm5ratperiods}
        \begin{split}
        &p^+(F_{k,\gamma}^+)(x) + p^-(F_{k,\gamma}^-)(x) \\
        &\quad \doteq -2\sum_{\substack{[a,b,c] \in [Q_\gamma] \\ a < 0 < c}}(ax^2-bx + c)^{k-1} - \frac{2D^{k-1/2}\zeta_{Q_\gamma}(k)}{\binom{2k-2}{k-1}(2k-1)\zeta(2k)}(x^{2k-2}-1),
                \end{split}
        \end{align}
where $\zeta_{Q_\gamma}(s)$ is the $\zeta$-function associated with $Q_\gamma$ as in \cite[p. 222]{kohnenzagier}, $\zeta(s)$ is the Riemann $\zeta$-function, and $\doteq$ means equality up to a non-zero multiplicative constant\footnote{The formulas \eqref{thm5ratperiods} and \eqref{thm4ratperiods} are correct if we normalize $F_{k,\gamma}$ as in \cite{kohnenzagier}. Since our normalization of $F_{k,\gamma}$ is different, we get some simple but unpleasant extra factors.}.
In particular, the periods $p_n(F_{k,\gamma}^+)$ for even $0< n < 2k-2$ and the periods $p_n(F_{k,\gamma}^-)$ for odd $0 < n < 2k-2$ are rational. 
 
We let $F_{k,D}(z) = \sum_{D_\gamma = D}F_{k,\gamma}(z)$ where the sum ranges over a system of representatives $\gamma$ of the conjugacy classes of primitive hyperbolic elements in $\SL$ with discriminant $D_\gamma = D$. From (\ref{thm5ratperiods}) it follows that for odd $k \in \Z$ we have  
\begin{align} \label{thm4ratperiods}
    \begin{split}
        p^+(F_{k, D})
        & \doteq -\sum_{\substack{[a,b,c] \in \mathcal{Q}_D \\ a < 0 < c}}(ax^2+bx + c)^{k-1} - \frac{D^{k-1/2}\zeta(k)L_D(k)}{\binom{2k-2}{k-1}(2k-1)\zeta(2k)}(x^{2k-2}-1),
        \end{split}
\end{align}
where $L_D(s)$ is the Dirichlet $L$-function associated to the Kronecker symbol $\left( \frac{D}{\cdot} \right)$. Formula (\ref{thm4ratperiods}) is the analog of \cite[Theorem 4]{kohnenzagier}. In particular, the periods $p_n(F_{k,D})$ for even $0 < n < 2k-2$ are rational and satisfy the symmetry $p_{2k-2-n}(F_{k, D}) = p_n(F_{k, D})$.

\subsection{The homogenized cycle integral as an "inner product"} By unfolding and \cite[Proposition~7]{kohnen}, one can also show that 
  \begin{align*}
   & \lim_{n \to \infty} \int_{\sigma^n.z_0}^{\sigma^{n + 1}.z_0} F_{k, \gamma}(z) Q_\sigma(z, 1)^{k - 1} dz \\
   &\qquad \doteq \int_{\G \setminus \uh} \sum_{g \in \Gamma_\gamma\backslash \Gamma} \sum_{h \in \Gamma_\sigma\backslash \Gamma} \frac{\sgn((Q_{\gamma}\circ g)(w_{h \sigma h^{-1}}', 1))}{(Q_\gamma \circ g)(z, 1)^k (Q_\sigma\circ h)(\overline{z}, 1)^k} y^{2k} \frac{dxdy}{y^2}.
    \end{align*}
    Since $w_{h \sigma h^{-1}}' = h^{-1}.w_{\sigma}'$, the integrand is $\G$-invariant. 
    
    One should compare this with the fact that the cycle integral $\displaystyle\int_{z_0}^{\sigma.z_0} f_{k, \gamma}(z) Q_\sigma(z, 1)^{k - 1} dz$ is up to constants equal to the Petersson inner product $\langle f_{k, \sigma}, f_{k, \gamma} \rangle$ of the hyperbolic Poincar\'e series (\ref{fkgamma}).

\subsection{Equidistribution of intersection angles} A possible application of Theorem \ref{mainresult} could be to give another proof the fact that the intersection angles $\theta_p$ of a fixed geodesic $[S_\gamma]$ with the geodesics of discriminant $D$ (for $\SL$) equidistribute to the measure $\frac{1}{2} \sin \theta \; d \theta$ as $D \to + \infty$. This was conjectured by Rickards \cite[Conjecture~4.2]{rickards} and recently proved by Jung and Sardari \cite{jungsardari}. To show this, it suffices to prove that $\cos \theta_p$ equidistribute to the Lebesgue measure on $[-1, 1]$. The Legendre polynomials $P_{k-1}$ with $k \geq 1$ form a complete orthonormal system of the space of continuous functions on $[-1, 1]$. For $k > 1$ even, the Weyl sum on the right hand side of (\ref{katokformula}) can be estimated using the Shimura-theory of cusp forms \cite{kohnen} and non-trivial bounds on the Fourier coefficients of half-integral weight cusp forms \cite{iwaniec}. A Siegel-type bound for the number of intersections can be obtained using (\ref{matsusakaformula}), \cite[Theorems~3.3 and 4.7, Remark~4.8]{matsusaka}, and some elementary considerations on continued fractions. For $k > 1$ odd, this approach does not work due to the appearance of the $\mu_p$-factor, so the Weyl sum here is given by the right hand side in (\ref{modintkatok}). However, it appears to be difficult to estimate the traces of the homogenized cycle integrals. We plan to come back to this in the future.

\subsection*{Acknowledgments} We are indebted to Özlem Imamo\={g}lu for suggesting the topic of the paper to us and for many insightful discussions. Moreover, we thank \'Arp\'ad T\'oth and Toshiki Matsusaka for helpful discussions. The first author was supported by SNF project 200021\_185014 and the the second author was supported by SNF projects 200021\_185014 and PZ00P2\_202210.

\end{document}